\documentclass{amsart}
\usepackage{xypic}
\xyoption{graph}

\newcommand{\QQ}{\mathbf{Q}}
\newcommand{\ZZ}{\mathbf{Z}}
\newcommand{\FF}{\mathbf{F}}
\newcommand{\OO}{\widetilde{O}}
\newcommand{\reduce}[1]{{\overline{#1}}}
\newcommand{\floor}[1]{\left\lfloor #1 \right\rfloor}
\newcommand{\abs}[1]{\left| #1 \right|}

\newcommand{\ep}{\varepsilon}
\newcommand{\error}[1]{\varepsilon({#1})}

\newcommand{\polytime}{{\mathsf M}}
\newcommand{\matrixtime}{{\mathsf{MM}}}

\newtheorem{thm}{Theorem}
\newtheorem{prop}[thm]{Proposition}
\newtheorem{lem}[thm]{Lemma}
\newtheorem{claim}{Claim}
\theoremstyle{definition}
\newtheorem{example}[thm]{Example}
\theoremstyle{remark}
\newtheorem*{rem}{Remark}

\begin{document}

\title{Kedlaya's Algorithm in Larger Characteristic}
\author{David Harvey}
\address{Department of Mathematics, Harvard University, 1 Oxford St, Cambridge MA 02138, USA}
\email{dmharvey@math.harvard.edu}
\urladdr{http://math.harvard.edu/\textasciitilde{}dmharvey/}

\begin{abstract}
We show that the linear dependence on $p$ of the running time of Kedlaya's point-counting algorithm in characteristic $p$ may be reduced to $p^{1/2}$.
\end{abstract}

\maketitle

\section{Introduction}
\label{sec:introduction}

In \cite{kedlaya}, Kedlaya introduced an algorithm for computing the zeta function of a hyperelliptic curve over $\FF_{p^n}$ of genus $g \geq 1$, which was remarkable for having running time polynomial in $g$ and $n$. Kedlaya did not discuss the dependence of the running time on $p$, and indeed at first it was thought that the algorithm would be practical only for very small primes. Later it was found that the dependence on $p$ was roughly linear (\cite{gaudry}, see also the survey paper \cite{kedlaya2}).

The main step of Kedlaya's algorithm --- the step where the linear dependence of the running time on $p$ occurs --- involves computing a $p$-adic approximation to the matrix of the $p$-th power Frobenius acting on a certain basis for the Monsky--Washnitzer cohomology of the curve (more precisely, the curve minus a few points). This is a $2g \times 2g$ matrix with entries in the degree $n$ unramified extension of $\QQ_p$. Kedlaya computes this matrix to precision $p^N$ in time $\OO(p N^2 g^2 n)$, where the `soft-oh' notation $\OO(X)$ indicates $O(X (\log X)^k)$ for some $k \geq 0$.

Our main result is the following. Let $\omega$ denote the exponent of matrix multiplication; that is, $\omega$ is a real number such that $m \times m$ matrices over a ring $R$ may be multiplied using $O(m^{\omega+\ep})$ ring operations in $R$ for any $\ep > 0$. Trivially one can take $\omega = 3$; see \cite{strassen} for the simplest example of a matrix multiplication algorithm that achieves $\omega < 3$.

\begin{thm}
\label{thm:complexity}
Let $N \geq 1$, and suppose that
\begin{equation}
\label{eq:p-assumption}
 p > (2N-1)(2g+1).
\end{equation}
Then the entries of the above matrix may be computed to precision $p^N$ in time
 \[ \OO(p^{1/2} N^{5/2} g^{\omega} n + N^4 g^4 n \log p). \]
In particular, for fixed $N$, $g$ and $n$, the running time is $\OO(p^{1/2})$.
\end{thm}

Our new algorithm is therefore superior to Kedlaya's original algorithm for fixed $g$ and $N$ and large enough $p$, but inferior for fixed $p$ and large enough $N$ or $g$. The final step of Kedlaya's algorithm is to compute the characteristic polynomial of the above matrix, but the running time of this step is only logarithmic in $p$, and will not concern us further.

The purpose of the assumption $p > (2N-1)(2g+1)$ is to simplify the analysis of denominators. It could be weakened somewhat, but the algorithm would become more complicated.

The author was motivated to develop this algorithm, not for point-counting purposes, but rather because of the role that the above matrix plays in the fast computation of $p$-adic heights of points on elliptic curves, as described in \cite{mst}. In that application, the parameter $N$ plays quite a different role. In \cite{kedlaya}, the aim is to compute the characteristic polynomial of Frobenius to sufficient precision that its exact value is pinned down by the Weil conjectures. Consequently Kedlaya takes $N = O(gn)$ and expresses all running time estimates in terms of $g$ and $n$ alone. On the other hand, in \cite{mst}, there is no reason to tie $N$ to $n$ or $g$. Indeed, $g = 1$ for an elliptic curve, and taking $n = 1$ suffices to handle curves defined over $\QQ$. Rather, the choice of $N$ ultimately depends on how accurately one wishes to determine the $p$-adic height. Therefore, in this paper we will analyse the dependence on $N$ separately from that of $n$ and $g$.

Our basic approach is the same as in \cite{kedlaya}: starting with an explicitly given basis of differentials for the Monsky--Washnitzer cohomology, we compute a representation of the action of an explicitly chosen lift of Frobenius on each basis differential, and then we apply a reduction algorithm that uses the cohomology relations to express the images as linear combinations of the original basis elements, thereby obtaining the desired matrix.

However, our algorithm differs from that of \cite{kedlaya} in two important respects. First, we make the key observation that the reductions in cohomology are given by formulae which may be interpreted as solving a \emph{linear recurrence with polynomial coefficients}. Therefore, instead of performing the reduction steps `one at a time', it becomes possible to use a baby-step/giant-step algorithm of Chudnovsky and Chudnovsky \cite{chudnovsky} to execute a whole \emph{sequence} of reductions in far less time than it would take to perform the reductions consecutively.

Second, to exploit this idea we must use a different representation for the relevant differentials. The difficulty is that in \cite{kedlaya}, the images of the basis differentials under Frobenius are approximated by series whose number of terms is at least linear in $p$, making it impossible to reach a running time proportional to $p^{1/2}$. We will use instead a different series approximation whose number of terms depends only on $g$ and $N$, not on $p$.

Rather than using the Chudnovskys' algorithm as they presented it, we take advantage of a modification due to Bostan, Gaudry and Schost \cite{bgs}, that improves on the running time by a factor logarithmic in the length of the recurrence. In our setting this translates to a speedup of $O(\log(pN))$, which for the feasible range of $p$ is very significant.

The relationship between our algorithm and the paper \cite{bgs} runs somewhat deeper. As one of the principal applications of their improved techniques for solving recurrences, they give an algorithm for computing the zeta function of a hyperelliptic curve over a finite field. Their approach is quite different to Kedlaya's, relying on the representation of the entries of the Hasse-Witt matrix associated to the curve $y^2 = f(x)$ as certain coefficients of the polynomial $f(x)^{(p-1)/2} \pmod p$. They then use the Chudnovskys' idea to efficiently compute those \emph{selected} coefficients, without computing the whole polynomial. It is striking that the Chudnovskys' algorithm plays such a central role in these two quite different approaches to computing zeta functions.

Our algorithm improves on the zeta function algorithm of \cite{bgs} in several ways, all of which may be traced to our essentially $p$-adic viewpoint. Whereas we obtain the zeta function modulo $p^N$ for any $N \geq 1$, their algorithm is only able to recover the zeta function modulo $p$, and they must then use other methods, such as $\ell$-adic methods, to obtain further information \cite[pp.~1800--1801]{bgs}. Furthermore, they achieve a running time of $\OO(p^{1/2} g^{3/2 + \omega} n)$ \cite[Theorem 17]{bgs}, which falls behind our algorithm by a factor of $g^{3/2}$ (ignoring the term involving $\log p$). The factor of $g^{3/2}$ may be accounted for as follows. In both our algorithm and the algorithm of \cite{bgs}, it is occasionally necessary to divide by $p$. To prevent precision loss at these division steps, \cite{bgs} are forced to lift from working modulo $p$ to working $p$-adically, artificially introducing $O(g)$ safety digits \cite[p.~1798]{bgs}. In our setting, the extra $p$-adic digits are ``already there'', and it is simply a matter of analysing the propagation of $p$-adic error terms. This explains a factor of $g$. The remaining factor of $g^{1/2}$ is more technical; essentially it occurs because our ``reduction matrices'' (see \S\ref{sec:reductions}) have certain $p$-adic analyticity properties that reduce the total number of matrices we must compute (see \S\ref{sec:compute-reduction-maps}).

Hubrechts \cite{hubrechts}, following a suggestion of Lauder, recently showed how to combine Kedlaya's algorithm with Dwork's deformation theory to improve the asymptotic running time with respect to $n$ (although the dependence on $g$ becomes worse). It would be interesting to study whether our approach to handling large $p$ is compatible with these developments.

\subsection*{Organisation of the paper}

In \S\ref{sec:notation} we fix notation, and in \S\ref{sec:kedlaya} we outline Kedlaya's original algorithm. In \S\ref{sec:frobenius-action} we give our alternative expression for the action of Frobenius on the appropriate differentials. In \S\ref{sec:reductions} we reformulate certain cohomological reductions as linear recurrences. In \S\ref{sec:recurrences} we give a slight generalisation of the algorithm of \cite{bgs} for solving linear recurrences. In \S\ref{sec:main} we describe the main algorithm, prove its correctness, and analyse its complexity. Finally, in \S\ref{sec:samples} we give some examples of timings for an implementation of the algorithm.

\subsection*{Acknowledgements}

Many thanks to Kiran Kedlaya for supplying the first clue that led to this algorithm, and for many helpful discussions about his algorithm, particularly regarding the thorny questions of precision loss. I would also like to thank William Stein for introducing me to the problem of computing $p$-adic heights, and for supplying the hardware on which the sample computations were performed (funded by NSF grant No.~0555776). Thanks to Barry Mazur, Kiran Kedlaya, Karim Belabas, William Stein, and an anonymous referee for several helpful comments on an early version of this paper.

\section{Notation and setup}
\label{sec:notation}

We will follow the notation of \cite{kedlaya} fairly closely. Let $p \geq 3$ be a prime, and let $q = p^n$ for some $n \geq 1$. The finite fields with $p$ and $q$ elements are denoted by $\FF_p$ and $\FF_q$. We denote by $\QQ_q$ the unramified extension of $\QQ_p$ of degree $n$, and by $\ZZ_q$ its ring of integers.

Let $\reduce Q \in \FF_q[x]$ be a monic polynomial of degree $2g+1$ ($g \geq 1$) with no multiple roots, so that the equation $y^2 = \reduce Q(x)$ defines the (projective) hyperelliptic curve $C/\FF_q$ of interest. We select an arbitrary lift $Q \in \ZZ_q[x]$ of $\reduce Q(x)$, also monic and of degree $2g+1$. (Note that in the application to computing $p$-adic heights \cite{mst}, the input data is actually $Q$ itself, rather than just $\reduce Q$.)

Let
 \[ \reduce A = \FF_q[x, y, y^{-1}]/(y^2 - \reduce Q(x)); \]
this is the coordinate ring of the curve $C'$ obtained from $C$ by removing the point at infinity and the points whose abscissae are the zeroes of $\reduce Q(x)$. Let
 \[ A = \ZZ_q[x, y, y^{-1}]/(y^2 - Q(x)) \]
be the lift of $\reduce A$ associated to $Q(x)$. Let $A^\dagger$ be the weak completion of $A$; explicitly, $A^\dagger$ is the ring of power series
 \[ \sum_{i \geq 0} \sum_{j \in \ZZ} a_{i, j} x^i y^j, \qquad a_{i, j} \in \ZZ_q, \]
such that $v_p(a_{i, j}) \to \infty$ at least linearly in $\max(i, \abs{j})$.

We will work mainly in the module $\Omega^-$ of differentials of $A^\dagger$ over $\QQ_q$ on which the hyperelliptic involution acts by $-1$. Explicitly, these are expressions of the form
 \[ \sum_{s \geq 0} \sum_{t \in \ZZ} a_{s, t} x^s y^{2t} dx/y,  \qquad a_{i, j} \in \QQ_q, \]
where the $a_{s, t}$ are subject to the same decay condition as above. Two differentials $\omega, \eta \in \Omega^-$ are cohomologous, denoted $\omega \sim \eta$, if there exists some $f \in A^\dagger \otimes \QQ_q$ such that $\omega - \eta = df$. We define the \emph{reduction} of $\omega$ to be the unique differential $\eta = B(x) dx/y$, cohomologous to $\omega$, such that the degree of $B \in \QQ_q[x]$ is at most $2g-1$. The existence and uniqueness of $\eta$ follows from the fact that $\{x^i dx/y\}_{i=0}^{2g-1}$ forms a basis for the Monsky--Washnitzer cohomology \cite[p.~329]{kedlaya}.

We lift the $p$-power Frobenius on $\FF_q$ to $A^\dagger$ as follows. On $\ZZ_q$, we take the canonical Witt vector Frobenius. We set $x^\sigma = x^p$,
\begin{equation}
\label{eq:frob-y}
 (y^{-1})^\sigma = y^{-p} \sum_{k=0}^\infty \binom{-1/2}{k} \frac{(Q(x)^\sigma - Q(x)^p)^k}{y^{2pk}},
\end{equation}
and $y^\sigma = (y^{-\sigma})^{-1}$. The above series converges in $A^\dagger$ (because $Q(x)^\sigma - Q(x)^p$ is divisible by $p$), and the definition ensures that $\sigma$ is an endomorphism of $A^\dagger$. We further extend $\sigma$ to $\Omega^-$ by $\sigma(f \, dg) = f^\sigma d(g^\sigma)$.

\section{A sketch of Kedlaya's original algorithm}
\label{sec:kedlaya}

In this section we will briefly describe Kedlaya's algorithm, paying particular attention to the dependence of the running time on $p$.

He begins by computing an approximation to $y^{-\sigma}$ of the form
 \[ y^{-\sigma} \approx y^{-p} \sum_{k=0}^{pN-1} \frac{A_k(x)}{y^{2k}}, \]
where each $A_k$ has degree at most $2g$. It is an approximation in two senses: it is truncated at a certain power of $y^{-2}$, and the coefficients are represented modulo $p^{N'}$, for some appropriately chosen $N'$ (slightly larger than $N$). Note that the time committed is already proportional to at least $p$, for the number of terms in the above series is about $Np$.

Next he takes the basis $\{x^i \, dx/y\}_{i=0}^{2g-1}$ for the de Rham cohomology of $A$ (actually, for its minus eigenspace under the hyperelliptic involution). Using the above series expansion of $y^{-\sigma}$, he computes an approximation to the image of each basis element under Frobenius,
\begin{equation}
\label{eq:frob-diff}
 \sigma(x^i \, dx/y) = x^{pi} d(x^p) y^{-\sigma} = p x^{pi+p-1} y^{-\sigma} dx
\end{equation}
as a series of the form
\begin{equation}
\label{eq:frob-diff-2}
 \sigma(x^i \, dx/y) \approx \sum_j \frac{F_j(x)}{y^{2j}} dx/y,
\end{equation}
where each $F_j$ has degree at most $2g$, and where again the series have about $Np$ terms.

For each $i$, he then applies a reduction algorithm to the terms on the right hand side of \eqref{eq:frob-diff-2}. At each step, he uses the identities $y^2 = Q(x)$ and $2y\, dy = Q'(x)dx$, together with the fact that $d(x^s y^t) = 0$ in cohomology for any $s$ and $t$, to reduce the term $F_j(x) y^{-2j} dx/y$ to a lower power of $y^{-2}$ (or in some cases, $y^2$). The terms are swept up sequentially until reaching $j = 0$. At this point one has computed the reduction of $\sigma(x^i dx/y)$, whose coefficients give the $(i+1)$-th column of the Frobenius matrix. The reduction step is performed once for each $j$, so again the total time is proportional to at least $p$.

\section{The Frobenius action on differentials}
\label{sec:frobenius-action}

As noted above, one of the barriers to making Kedlaya's algorithm run in time less than linear in $p$ is that the series approximation for $\sigma(x^i dx/y)$ given by \eqref{eq:frob-diff-2} has about $Np$ terms. The following proposition gives a different approximation for $\sigma(x^i dx/y)$ that requires only $O(N^2 g)$ terms; in particular, the number of terms does not depend on $p$.

\begin{prop}
\label{prop:frob-diff-approx}
Let $C_{j, r} \in \ZZ_q$ be the coefficient of $x^r$ in $Q(x)^j$. For $0 \leq j < N$, let
 \[ B_{j, r} = p C_{j, r}^\sigma \sum_{k=j}^{N-1} (-1)^{k+j} \binom{-1/2}{k} \binom{k}{j} \quad \in \ZZ_q. \]
For $0 \leq i < 2g$, set
\begin{equation}
\label{eq:frob-diff-approx}
 T_i = \sum_{j=0}^{N-1} \sum_{r=0}^{(2g+1)j} B_{j, r} x^{p(i+r+1)-1} y^{-p(2j+1)+1} dx/y.
\end{equation}
Then the reduction of $T_i$ agrees modulo $p^N$ with the reduction of $\sigma(x^i dx/y)$.
\end{prop}

\begin{proof}
From \eqref{eq:frob-y} and \eqref{eq:frob-diff} we obtain
\begin{equation}
\label{eq:frob-diff-expand}
 \sigma(x^i dx/y) = \sum_{k=0}^\infty p \binom{-1/2}{k} (Q(x)^\sigma - Q(x)^p)^k x^{pi+p-1} y^{-p(2k+1)+1} dx/y.
\end{equation}
Since $Q(x)^\sigma - Q(x)^p$ is divisible by $p$, the $k$-th term $U_k$ of \eqref{eq:frob-diff-expand} is of the form
 \[ U_k = p^{k+1} F(x) y^{-p(2k+1)} dx, \]
where $F \in \ZZ_q[x]$ has degree at most
 \[ ((2g+1)p-1)k + pi + p - 1 < (2g+1)(k+1)p. \]
By repeatedly dividing $F(x)$ by $Q(x) = y^2$, we may rewrite this as
 \[ U_k = p^{k+1} \sum_{j=0}^{(k+1)p - 1} F_j(x) y^{-p(2k+1)+2j} dx, \]
where each $F_j \in \ZZ_q[x]$ has degree at most $2g$.

We must show that the coefficients of the reduction of $U_k$ are divisible by $p^N$, for all $k \geq N$. The terms for which $0 \leq j < (k+\frac12)p$ may be handled by \cite[Lemma 2]{kedlaya}, which shows that the reduction of $F_j(x) y^{-p(2k+1)+2j}$ becomes integral on multiplication by $p^{1+\lfloor \log_p(2k+1)\rfloor}$. Assumption \eqref{eq:p-assumption} implies that $\floor{\log_p(2k+1)} \leq k - N$, which covers this case. The remaining terms for which $(k+\frac12)p \leq j \leq (k+1)p - 1$ require \cite[Lemma 3]{kedlaya}. (Note: Lemma 3 as stated in \cite{kedlaya} is incorrect. A corrected version is in the errata to \cite{kedlaya}, and a proof is given in Lemma 4.3.5 of \cite{edixhoven}.) For these $j$ we find that the reduction of $F_j(x) y^{-p(2k+1)+2j} dx$ becomes integral on multiplication by $p^m$ where
 \[ m = \floor{\log_p((2g+1)(-p(2k+1)+2j+2)-2)} \leq \floor{\log_p((2g+1)p)} \leq 1, \]
the last inequality again depending on \eqref{eq:p-assumption}.

Consequently the terms in \eqref{eq:frob-diff-expand} for $k \geq N$ do not contribute modulo $p^N$ to the reduction of $\sigma(x^i dx/y)$, so we may ignore them. Therefore, let
 \[ T_i = \sum_{k=0}^{N-1} p \binom{-1/2}{k} (Q(x)^\sigma - Q(x)^p)^k x^{pi+p-1} y^{-p(2k+1)+1} dx/y. \]
We now replace $Q(x)^p$ by $y^{2p}$, use the binomial formula to expand $(Q(x)^\sigma - y^{2p})^k$, and write out the coefficients $Q(x)^\sigma$ explicitly in terms of the $C_{j, r}$. After rearranging the summations, we obtain the representation for $T_i$ indicated in the statement of the proposition.
\end{proof}

\begin{rem}
Ultimately, the linear contribution of $p$ to the running time of Kedlaya's original algorithm arises from explicitly expanding out the $Q(x)^p$ term in a formula of the above type. In the proof of Proposition \ref{prop:frob-diff-approx}, we avoided this by substituting $y^{2p}$ for $Q(x)^p$, and we will see that our algorithm will accordingly never need to compute the coefficients of $Q(x)^p$. At first glance this may seem odd, since in Kedlaya's original algorithm, the expansion of $Q(x)^p$ --- more precisely, the congruence modulo $p$ between $Q(x)^p$ and $Q(x)^\sigma$ --- is precisely what causes the terms in $y^\sigma$ with high powers of $y^{-2}$ to have $p$-adically small coefficients. In our case however, one finds that the reduction of each term $B_{j, r} x^{p(i+r+1)-1} y^{-p(2j+1)+1} dx/y$ of $T_i$ generally contributes to \emph{all $N$ digits} of the coefficients of the reduction of $T_i$, regardless of the value of $r$ or $j$. In fact, even the \emph{sum} of all terms for a given power of $y^{-2}$ (that is, for a given $j$) contributes to all $N$ digits. It is almost as if our algorithm ignores the decay conditions defining $A^\dagger$. Of course those decay conditions do play a role, by inducing hidden cancellations among the $B_{j, r}$.
\end{rem}

\section{Horizontal and vertical reduction}
\label{sec:reductions}

Let $s \geq -1$ and $t \in \ZZ$. We define $W_{s, t}$ to be the $\QQ_q$-vector space of differentials of the form
 \[ F(x) x^s y^{-2t} dx/y, \]
where $F(x) \in \QQ_q[x]$ has degree at most $2g$. In the case $s = -1$, we impose the additional condition that the constant term of $F(x)$ must be zero (so that none of the differentials ever involve negative powers of $x$).

In \S\ref{sec:vertical} and \S\ref{sec:horizontal} we will give maps between the various $W_{s, t}$ that send differentials to cohomologous differentials. The point is to give explicit formulae, so that the maps may be interpreted as defining linear recurrences. First we discuss `vertical' reductions, which map $W_{-1, t}$ to $W_{-1, t-1}$; this is the main type of reduction that appears in \cite{kedlaya}. Then we discuss `horizontal' reductions, which map $W_{s, t}$ to $W_{s-1, t}$. The aim is to eventually reduce everything to $W_{-1, 0}$, since this space consists of the differentials of the form $G(x) dx/y$, where $G$ has degree at most $2g-1$.

We will generally identify elements of $W_{s, t}$ with vectors in $\ZZ_q^{2g+1}$ (or $\ZZ_q^{2g}$ in the case $s = -1$), with respect to the basis $\{x^{i+s} y^{-2t} dx/y\}_{i=0}^{2g}$ (or with respect to $\{x^i y^{-2t} dx/y\}_{i=0}^{2g-1}$ in the case $s = -1$).

\subsection{Vertical reduction}
\label{sec:vertical}

Let $0 \leq i < 2g$ and $t \in \ZZ$. Since $\reduce Q(x)$ has no repeated roots, we can find polynomials $R_i, S_i \in \ZZ_q[x]$, where $\deg R_i \leq 2g-1$ and $\deg S_i \leq 2g$, such that
\begin{equation}
\label{eq:bezout}
 x^i = R_i(x)Q(x) + S_i(x)Q'(x).
\end{equation}
(To get the integrality of $R_i$ and $S_i$, we have used the assumption that $p > 2g+1$, so that the leading coefficient of $Q'(x)$ is a unit.) Using the relation $2y \, dy = Q'(x) dx$, we have
 \[ x^i y^{-2t} dx/y = R_i(x) y^{-2t+2} dx/y + 2 S_i(x) y^{-2t} dy. \]
Since $d(S_i(x) y^{-2t+1})$ is zero in cohomology, after a little algebra we find that
\begin{equation}
\label{eq:vertical-relation}
 x^i y^{-2t} dx/y \sim \frac{(2t-1)R_i(x) + 2S_i'(x)}{2t-1} y^{-2t+2} dx/y.
\end{equation}
(The above calculation is essentially the one in \cite[p.~329]{kedlaya}.)

This last relation may be rephrased in terms of the vector spaces $W_{-1, t}$ as follows.

\begin{prop}
\label{prop:single-vertical-reduction}
Let
 \[ M_V(t) : W_{-1, t} \to W_{-1, t-1} \]
be the linear map given by the $2g \times 2g$ matrix whose $(i+1)$-th column consists of the coefficients of the polynomial $(2t-1)R_i(x) + 2S_i'(x)$. Let
 \[ D_V(t) = 2t - 1. \]
Then for any $\omega \in W_{-1, t}$, we have
 \[ \omega \sim D_V(t)^{-1} M_V(t) \, \omega \quad (\in W_{-1, t-1}). \]
\end{prop}
In other words, $D_V(t)^{-1} M_V(t)$ is the \emph{reduction matrix} for transporting a differential from $W_{-1, t}$ to a cohomologous differential in $W_{-1, t-1}$. Note that the entries of $M_V(t)$ are linear polynomials in $\ZZ_q[t]$, as is $D_V(t)$.

We will be interested in iterating this process. For $t_0 < t_1$, let
 \[ M_V(t_0, t_1) : W_{-1, t_1} \to W_{-1, t_0} \]
be defined by
 \[ M_V(t_0, t_1) = M_V(t_0+1) M_V(t_0+2) \cdots M_V(t_1). \]
Similarly let
 \[ D_V(t_0, t_1) = D_V(t_0+1) D_V(t_0+2) \cdots D_V(t_1). \]
With this notation we obtain:
\begin{prop}
\label{prop:vertical-reduction}
For any $\omega \in W_{-1, t_1}$,
 \[ \omega \sim D_V(t_0, t_1)^{-1} M_V(t_0, t_1) \, \omega \quad (\in W_{-1, t_0}). \]
\end{prop}

\begin{example}[An elliptic curve]
\label{ex:elliptic}

We compute $M_V(t)$ for the elliptic curve $y^2 = Q(x) = x^3 + ax + b$. First solve \eqref{eq:bezout} for $i = 0, 1$, obtaining
\[
\begin{aligned}
 R_0(x) & = \Delta^{-1} (-18ax + 27b) \\
 S_0(x) & = \Delta^{-1} (6ax^2 - 9bx + 4a^2) \\
 R_1(x) & = \Delta^{-1} (27bx + 6a^2) \\
 S_1(x) & = \Delta^{-1} (-9bx^2 - 2a^2x - 6ab),
\end{aligned}
\]
where $\Delta = 27b^2 + 4a^3$ is the discriminant of the curve. Therefore
\[
\begin{aligned}
 (2t-1)R_0(x) + 2 S_0'(x)
    & = \Delta^{-1} (-6ax(6t - 7) + 9b(6t - 5))   \\
 (2t-1)R_1(x) + 2 S_1'(x)
    & = \Delta^{-1} (9bx(6t - 7) + 2a^2(6t - 5)),   \\
\end{aligned}
\]
and so the matrix $M_V(t)$ is given by
\[ M_V(t) = \Delta^{-1}
\begin{pmatrix}
 9b(6t-5)  &  2a^2(6t-5)  \\
 -6a(6t-7)  &  9b(6t-7)
\end{pmatrix}. \]
\end{example}

\subsection{Horizontal reduction}
\label{sec:horizontal}

Let $s \geq 0$ and $t \in \ZZ$. In cohomology,
\[
\begin{aligned}
0 & \sim d(x^s y^{-2t+1}) \\
  & = sx^{s-1} y^{-2t+1} dx - (2t - 1)x^s y^{-2t} dy \\
  & = \left( sx^{s-1} Q(x) - \frac12 (2t - 1) x^s Q'(x) \right) y^{-2t} dx/y.
\end{aligned}
\]
Decompose $Q(x)$ as
 \[ Q(x) = x^{2g+1} + P(x), \]
where $P \in \ZZ_q[x]$ has degree at most $2g$. After substituting this into the previous equation and rearranging, we obtain
\begin{equation}
\label{eq:horizontal-relation}
 x^{s + 2g} y^{-2t} dx/y \sim
    \frac{2sP(x) - (2t-1)x P'(x)}{(2g+1)(2t-1) - 2s}
        x^{s-1} y^{-2t} dx/y.
\end{equation}

\begin{prop}
\label{prop:single-horizontal-reduction}
Let
 \[ M_H^t(s) : W_{s, t} \to W_{s-1, t} \]
be the linear map given by the matrix
 \[
  M_H^t(s) = \begin{pmatrix}
   0         & 0         & \cdots & 0         & C_0(s) \\
   D_H^t(s)  & 0         &        & 0         & C_1(s) \\
   0         & D_H^t(s)  &        & 0         & C_2(s) \\
   \vdots    &           & \ddots &           & \vdots \\
   0         & 0         & \cdots & D_H^t(s)  & C_{2g}(s)
    \end{pmatrix},
 \]
where
 \[ D_H^t(s) = (2g+1)(2t-1) - 2s, \]
and where $C_h(s)$ is the coefficient of $x^h$ in the polynomial
 \[ C(x, s) = 2sP(x) - (2t-1)x P'(x). \]
Then for any $\omega \in W_{s, t}$, we have
 \[ \omega \sim D_H^t(s)^{-1} M_H^t(s) \, \omega \quad (\in W_{s-1, t}). \]
\end{prop}
\begin{proof}
The bulk of the statement follows from \eqref{eq:horizontal-relation}. In addition, the constant term of $C(x, 0)$ is zero, so $M_H^t(0)$ does indeed map into $W_{-1, t}$.
\end{proof}

Note that, for a fixed choice of $t$, the entries of $M_H^t(s)$ and $D_H^t(s)$ are linear polynomials in $\ZZ_q[s]$, and $D_H^t(s)$ does not vanish for any $s$, since it is always odd.

To iterate this process, we define, for $-1 \leq s_0 < s_1$,
 \[ M_H^t(s_0, s_1) : W_{s_1, t} \to W_{s_0, t} \]
by
 \[ M_H^t(s_0, s_1) = M_H^t(s_0+1) M_H^t(s_0+2) \cdots M_H^t(s_1), \]
and
 \[ D_H^t(s_0, s_1) = D_H^t(s_0+1) D_H^t(s_0+2) \cdots D_H^t(s_1). \]
We obtain:
\begin{prop}
\label{prop:horizontal-reduction}
For any $\omega \in W_{s_1, t}$,
 \[ \omega \sim D_H^t(s_0, s_1)^{-1} M_H^t(s_0, s_1) \, \omega \quad (\in W_{s_0, t}). \]
\end{prop}

\begin{example}[An elliptic curve]
We compute $D_H^t(s)$ and $M_H^t(s)$ for the elliptic curve $y^2 = Q(x) = x^3 + ax + b$. We have
 \[ D_H^t(s) =  6t - 2s - 3, \]
and $P(x) = ax + b$, so
 \[ 2sP(x) - (2t-1)xP'(x) = ax(2s - 2t + 1) + bs. \]
Then $M_H^t(s)$ is given by
 \[
  M_H^t(s) = \begin{pmatrix}
   0            & 0            & 2bs \\
   6t - 2s - 3  & 0            & a(2s - 2t + 1) \\
   0            & 6t - 2s - 3  & 0
    \end{pmatrix}.
 \]
\end{example}

\section{Algorithms for linear recurrences}
\label{sec:recurrences}

The following theorem from \cite{bgs} is not precisely what we will need, but it is close enough that we will be able to adapt it without difficulty. To state it, we need to introduce some notation from \cite{bgs}. Let $R$ be a commutative ring with identity. In this section we will work in an algebraic model of computation, so running times are measured by counting ring operations in $R$. We denote by $\polytime(d)$ the time required to multiply polynomials of degree $d$ over $R$, and by $\matrixtime(m)$ the time required to multiply $m \times m$ matrices with entries in $R$. In \cite{bgs} they make several reasonable regularity assumptions about the growth of $\polytime(d)$ and $\matrixtime(m)$, which are certainly satisfied in the cases we will consider.

For any integer $s \geq 0$, they define a certain quantity $\mathsf{D}(1, 2^s, 2^s) \in R$. The definition is straightforward, but lengthy, and we will not give it here. The only fact we need (see \cite[p.~1787]{bgs}) is that if $2, 3, \ldots, 2^s + 1$ are units in $R$, then $\mathsf{D}(1, 2^s, 2^s)$ is invertible in $R$, and that its inverse may be used to efficiently recover the inverses of certain other elements of $R$ that are needed in the interpolation steps of their algorithm.

Now, let $M(X)$ be an $m \times m$ matrix of linear polynomials in $R[X]$, and let $K \geq 1$ be an integer. Given an initial vector $U_0 \in R^m$, they define a sequence of vectors by
 \[ U_{i+1} = M(i+1) U_i \]
for $i \geq 0$. If one wishes to compute several $U_i$ in the range $0 \leq i \leq K$, the naive algorithm (simply iterating the above relation) requires time $O(m^2 K)$. The following theorem improves substantially on this, as long as not too many $U_i$ are requested.

\begin{thm}[\protect{\cite[Theorem 15]{bgs}}]
\label{thm:bgs}
Let $0 < K_1 < K_2 < \cdots < K_r = K$ be integers, and let $s = \floor{\log_4 K}$. Suppose that $2, 3, \ldots, 2^s + 1$ are invertible in $R$, and that the inverse of $\mathsf{D}(1, 2^s, 2^s)$ is known. Suppose also that $r < K^{\frac12 - \ep}$, with $0 < \ep < 1/2$. Then $U_{K_1}, \ldots, U_{K_r}$ can be computed using
 \[ O(\matrixtime(m) \sqrt K + m^2 \polytime(\sqrt K)) \]
ring operations in $R$.
\end{thm}

The theorem we require is a little stronger. Using similar notation to the horizontal and vertical reduction matrices of \S\ref{sec:reductions}, we define
 \[ M(k, k') = M(k') M(k'-1) \cdots M(k+2) M(k+1) \]
for $k < k'$. (Note that we have switched the ordering of the matrices from \S\ref{sec:reductions}, to match the notation of \cite{bgs}. It is trivial to adapt the algorithm to work in the opposite direction.) Instead of just computing the images $U_{K_1}, \ldots, U_{K_r}$ of a single vector $U_K$, our aim is to compute the matrices $M(K_i, L_i)$ for a sequence of intervals $(K_i, L_i)$. The following slight generalisation of Theorem \ref{thm:bgs} achieves this.

\begin{thm}
\label{thm:bgs-generalisation}
Let
 \[ 0 \leq K_1 < L_1 \leq K_2 < L_2 \leq \cdots \leq K_r < L_r \leq K \]
be integers, and let $s = \floor{\log_4 K}$. Suppose that $2, 3, \ldots, 2^s + 1$ are invertible in $R$, and that the inverse of $\mathsf{D}(1, 2^s, 2^s)$ is known. Suppose also that $r < K^{\frac12 - \ep}$, with $0 < \ep < 1/2$. Then $M(K_1, L_1), \ldots, M(K_r, L_r)$ can be computed using
 \[ O(\matrixtime(m) \sqrt K + m^2 \polytime(\sqrt K)) \]
ring operations in $R$.
\end{thm}
\begin{rem}
When we prove the main complexity result (Theorem \ref{thm:complexity}) we will ignore the distinction between the two terms in the above estimate. The key point is that the running time is soft-linear in $\sqrt K$, and polynomial in $m$.
\end{rem}
\begin{proof}
The algorithm is almost exactly the same as the one given in the proof of \cite[Theorem 15]{bgs}, so we will not spell out all the details. To explain it, we first give a very high-level sketch of their algorithm. In ``Step 0'', they compute a sequence of matrices
\begin{equation}
\label{eq:bgs-sequence}
 M(0, H), M(H, 2H), \ldots, M((B-1)H, BH),
\end{equation}
where both $H$ and $B$ are a small constant factor away from $\sqrt K$. They apply these matrices successively to $U_0$ to compute $U_{kH}$ for all $0 \leq k \leq B$. Each target index $K_i$ will fall within one of the intervals $[kH, (k+1)H]$. Then they perform a ``refining'' step, where they deduce $U_{K_i}$ from $U_{kH}$ by evaluating appropriate products of $M(X)$ over (much smaller) subintervals of $[kH, K_i]$. To stay within the time bounds, they use multipoint evaluation techniques to refine towards all target indices simultaneously.

(The main reason that their algorithm is a logarithmic factor faster than the Chudnovskys' algorithm is that in Step 0, they give up some control over which intervals are computed, in exchange for having available a faster method for computing them. This is why the separate refining step is necessary.)

To adapt this to our needs, we need only perform a little extra work. Given the input indices $K_i$ and $L_i$, we compute the sequence \eqref{eq:bgs-sequence}, using the same method as \cite{bgs}. We now perform a refining step using the same algorithm as in \cite{bgs}, but we will need to refine over more intervals. Suppose that $K_i$ lies in $[k_1 H, (k_1 + 1) H]$ and that $L_i$ lies in $[k_2 H, (k_2 + 1) H]$, where $k_2 \geq k_1$. If $k_1 = k_2$ then we refine over $[K_i, L_i]$. If $k_2 > k_1$, we must refine over both $[K_i, (k_1+1)H]$ and $[k_2 H, L_i]$.

After computing the products $M(k, k')$ for each of these intervals, we must perform an additional `gluing' step. Namely, each of our target intervals $(K_i, L_i)$ is a union of intervals $(k, k')$ for which $M(k, k')$ has been computed (either in Step 0 or in the refining step), and so we simply multiply together the $M(k, k')$ for those intervals, in the appropriate order.

To estimate the total time, we note first that our `Step 0' is identical to their `Step 0'. The refining steps take at most twice as long as theirs, since we have at most doubled the number of intervals to be considered, and the lengths of those intervals satisfy the same bounds. One must also check the invertibility conditions in $R$; these are still satisfied since they depend only on the maximum length of the intervals, which has not changed. Finally, the extra gluing step consists of at most $O(\sqrt K)$ matrix multiplications, costing time $O(\matrixtime(m) \sqrt K)$, which fits within the required time bound.
\end{proof}

\section{The main algorithm}
\label{sec:main}

In this section we describe the main algorithm for computing the Frobenius matrix. The basic idea is to start with the approximation $T_i$ for $\sigma(x^i dx/y)$ given by Proposition \ref{prop:frob-diff-approx}, and then to use the reduction maps to push each term towards $W_{-1, 0}$. Theorem \ref{thm:bgs-generalisation} is used to efficiently compute the reduction maps.

Figure \ref{fig:reductions} illustrates the strategy in the case $g = 1$ and $N = 3$. Each vertex corresponds to a $W_{s, t}$, and the arrows correspond to horizontal and vertical reductions. The black vertices are those which are the starting point for at least one term from some $T_i$. (There are additional vertices and arrows used in the algorithm that for reasons of clarity are not shown on the diagram.)

\newcommand{\white}{\circ}
\newcommand{\black}{\bullet}
\begin{figure}[h]
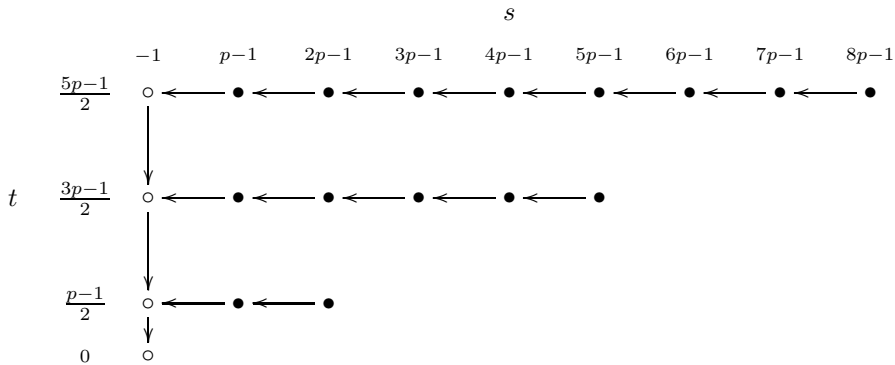

\xygraph{
!{<0cm,0cm>;<1.2cm,0cm>:<0cm,0.7cm>::} 
!{(0,5) }*+{\white}="05"
!{(1,5) }*+{\black}="15"
!{(2,5) }*+{\black}="25"
!{(3,5) }*+{\black}="35"
!{(4,5) }*+{\black}="45"
!{(5,5) }*+{\black}="55"
!{(6,5) }*+{\black}="65"
!{(7,5) }*+{\black}="75"
!{(8,5) }*+{\black}="85"
!{(0,3) }*+{\white}="03"
!{(1,3) }*+{\black}="13"
!{(2,3) }*+{\black}="23"
!{(3,3) }*+{\black}="33"
!{(4,3) }*+{\black}="43"
!{(5,3) }*+{\black}="53"
!{(0,1) }*+{\white}="01"
!{(1,1) }*+{\black}="11"
!{(2,1) }*+{\black}="21"
!{(0,0) }*+{\white}="00"
!{(-0.7,5)}*+{\frac{5p-1}2}
!{(-0.7,3)}*+{\frac{3p-1}2}
!{(-0.7,1)}*+{\frac{p-1}2}
!{(-0.7,0)}*+{\scriptstyle 0}
!{(0, 5.7)}*+{\scriptstyle -1}
!{(1, 5.7)}*+{\scriptstyle p-1}
!{(2, 5.7)}*+{\scriptstyle 2p-1}
!{(3, 5.7)}*+{\scriptstyle 3p-1}
!{(4, 5.7)}*+{\scriptstyle 4p-1}
!{(5, 5.7)}*+{\scriptstyle 5p-1}
!{(6, 5.7)}*+{\scriptstyle 6p-1}
!{(7, 5.7)}*+{\scriptstyle 7p-1}
!{(8, 5.7)}*+{\scriptstyle 8p-1}
!{(-1.5, 3)}*+{t}
!{(4, 6.5)}*+{s}
"15":"05" "25":"15" "35":"25" "45":"35"
"55":"45" "65":"55" "75":"65" "85":"75"
"13":"03" "23":"13" "33":"23" "43":"33" "53":"43"
"11":"01" "21":"11"
"01":"00" "03":"01" "05":"03"
}
\caption{Reduction strategy for $g = 1$ and $N = 3$}
\label{fig:reductions}
\end{figure}

One of the more magical aspects of Kedlaya's original algorithm is the way that $p$-adic precision losses propagate through the calculation. Although one needs to perform about $N$ divisions by $p$, Kedlaya shows that in fact only $O(\log_p N)$ spare digits of precision must be carried.

A similar argument applies to our algorithm, and since we have assumed $p$ to be sufficiently large compared to $g$ and $N$, it turns out that only one spare digit is necessary. However, some caution is required. For example, the product of all the $M_H^t(s)$ across a whole `row' of the horizontal reductions will generally be \emph{zero} modulo $p^N$, and therefore one must interleave the multiplications by $M_H^t(s)$ and divisions by $D_H^t(s)$ in such a way that the denominators can ``catch up with'' the build-up of $p$-divisibility of the numerators. In \S\ref{sec:horizontal-algorithm} we perform a more detailed analysis, showing how to do almost all of the work with \emph{no} spare digits at all. In practical terms, avoiding even this single extra digit yields enormous savings in time and memory when $N$ is small. For the vertical reductions, at least in the case $N > 1$, this kind of analysis seems much more difficult, and consequently we will retain the spare digit.

\subsection{Preliminaries}

The algorithm works in two different rings, $R_0 = \ZZ_q/(p^N)$ and $R_1 = \ZZ_q/(p^{N+1})$. At certain stages we will need to compute $a/b$, where $b$ is not a unit; we may take the result to be any $c$ satisfying $bc = a$. We will see below that such divisions will always be possible in $\ZZ_q$ when they occur, and that the errors introduced do not contribute to the final result modulo $p^N$.

As a preliminary step, we compute the coefficients $B_{j, r}$ given in Proposition \ref{prop:frob-diff-approx}, for $0 \leq j < N$ and $0 \leq r \leq (2g+1)j$, as elements of $R_1$.

Let us write $T_i$ as
\[
\begin{aligned}
 T_i & = \sum_{j=0}^{N-1} T_{i, j}, \qquad T_{i, j} = \sum_{k=0}^{i+(2g+1)j+1} T_{i, j, k}, \\
 T_{i, j, k} & = B_{j, k-i-1} x^{pk-1} y^{-p(2j+1)+1} dx/y,
\end{aligned}
\]
where for convenience we declare that $B_{j, r} = 0$ for $r < 0$. Note that $T_{i, j, k} \in W_{pk-1, t}$, where $t = \frac12((2j+1)p - 1)$.

\subsection{Horizontal reduction phase}
\label{sec:horizontal-algorithm}

This phase is performed once for each $0 \leq j < N$; throughout this section we regard $j$ as fixed.

Let $t = \frac12((2j+1)p - 1)$. The aim is to use the horizontal reduction maps to find differentials $w_{i, j} \in W_{-1, t}$ that are cohomologous to $T_{i, j}$, and whose coefficients are correct modulo $p^N$, for $0 \leq i < 2g$.

\subsubsection{Computing the reduction maps}
\label{sec:compute-reduction-maps}

Let $L = (2g+1)j + 2g$. We must first compute the horizontal reduction matrices
\begin{equation}
\label{eq:horizontal-blocks}
\begin{aligned}
  M(k) & = M^t_H((k-1)p, kp-2g-2), \\
  D(k) & = D^t_H((k-1)p, kp-2g-2),
\end{aligned}
\end{equation}
for $1 \leq k \leq L$, with entries in $R_0$. (Once computed, it may be convenient to lift them to $R_1$, but it is only necessary to know them modulo $p^N$.)

This is accomplished in two steps. We will discuss $M(k)$ only; the $D(k)$ are handled entirely analogously.

The first and most time-consuming step is to use Theorem \ref{thm:bgs-generalisation} to compute $M(k)$ for $1 \leq k \leq L'$, where $L' = \min(N, L)$. To verify the invertibility hypotheses of Theorem \ref{thm:bgs-generalisation}, we must check that $\sqrt K + 1 < p$, where $K = L'p - 2g - 2$ is the total length of the interval containing all the reduction intervals. From \eqref{eq:p-assumption} we know that $(2g+1)(2N-1) \leq p-1$, so
\[
\begin{aligned}
 2K & \leq (2g+1)(2N-1)p + (2g+1)p - 2p - 4g - 4 \\
    & \leq (p-1)p + (p-1)p - 2p - 4g - 4 \\
    & < 2(p-1)^2,
\end{aligned}
\]
from which the desired inequality follows.

The second step is to deduce the remaining $M(k)$ for $N < k \leq L$. (This is of course only necessary when $L > N$.) It is possible to simply use Theorem \ref{thm:bgs-generalisation} again, but it is much more efficient to take advantage of the known values $M(1), \ldots M(N)$. If $N = 1$ this is trivial, since the $M(k)$ are all equal modulo $p$. The author thanks Kiran Kedlaya for suggesting the following interpolation method to handle the case $N > 1$.

Consider the matrix
 \[ F(s) = M^t_H(s - p + 1) \cdots M^t_H(s - 2g - 2), \]
which is a matrix of polynomials in $s$. Expanding as a Taylor series in $s$, we obtain
 \[ M(k) = F(kp) = F(0) + F'(0) kp + \cdots + \frac{1}{(N-1)!} F^{(N-1)}(0) (kp)^{N-1} \pmod{p^N}. \]
Then by simple linear algebra, the values of $F(kp) \pmod{p^N}$ for $1 \leq k \leq N$ determine completely the values of $F^{(i)}(0) p^i / i!$ for $0 \leq i < N$. Namely, we have
 \[
  \begin{pmatrix} F(p) \\ F(2p) \\ \vdots \\ F(Np) \end{pmatrix}
  =
  \begin{pmatrix}
  1      & 1  & \cdots & 1 \\
  1      & 2  &        & 2^{N-1} \\
  \vdots & &  & \vdots \\
  1      & N  & \cdots & N^{N-1}
  \end{pmatrix}
  \begin{pmatrix} F(0) \\ F'(0)p \\ \vdots \\ \frac1{(N-1)!} F^{(N-1)}(0)p^{N-1}, \end{pmatrix}
 \]
and the Vandermonde matrix is invertible modulo $p$ (since $p > N$). After solving for the $F^{(i)}(0) p^i / i!$, the remaining $M(k)$ are computed by substituting the appropriate values of $k$ into the above Taylor series.

\begin{rem}
In the case $N = 1$ there is a yet faster method available for computing $D(k)$ (although not $M(k)$). Namely, since $t \equiv -1/2 \pmod p$ we have
 \[ D(k) \equiv D(1) \equiv \prod_{s=1}^{p-2g-2} -2(2g+1) -2s \pmod p, \]
which by Wilson's theorem is equal to
 \[ (-2)^{p - 2g - 2} \prod_{s=2g+2}^{p-1} s \equiv (2^{2g+1} (2g+1)!)^{-1} \pmod p. \]
\end{rem}

\subsubsection{Performing the reductions}

Now we fix $0 \leq i < 2g$, and show how to compute $w_{i, j}$. We will define a sequence of differentials $v_{i+(2g+1)j+1}, \ldots, v_1, v_0$, where $v_m \in W_{mp-1, t}$, with the property that
\begin{equation}
\label{eq:horizontal-pile}
 v_m \sim \sum_{k=m}^{i+(2g+1)j+1} T_{i, j, k}.
\end{equation}
In particular we will have $v_0 \sim T_{i, j}$, so this $v_0$ is the $w_{i, j}$ that we seek. The $v_m$ are computed with entries in $R_1$. However, not all their $p$-adic digits will be correct; we will say more about this in a moment.

Naturally, the sequence begins with
 \[ v_{i+(2g+1)j+1} = T_{i, j, i+(2g+1)j+1}. \]
Then, given $v_m$, we compute $v_{m-1}$ as follows. We first move from $W_{mp-1,t}$ to $W_{mp-2g-2,t}$, one step at a time, via the following sequence:
 \[ \begin{aligned}
 v_m^{(1)} & = v_m & & \in W_{mp-1,t}, \\
 v_m^{(2)} & = D^t_H(mp-1)^{-1} M^t_H(mp-1) v_m^{(1)} & & \in W_{mp-2, t}, \\
 v_m^{(3)} & = D^t_H(mp-2)^{-1} M^t_H(mp-2) v_m^{(2)} & & \in W_{mp-3, t}, \\
 & \vdots \\
 v_m^{(2g+2)} & = D^t_H(mp-2g-1)^{-1} M^t(mp-2g-1) v_m^{(2g+1)} & & \in W_{mp-2g-2, t}.
 \end{aligned} \]
Using the reduction matrices \eqref{eq:horizontal-blocks} computed above, we set
\begin{equation}
\label{eq:transport-to-v-prime}
 v_m' = D^t_H((m-1)p, mp-2g-2)^{-1} M^t_H((m-1)p, mp-2g-2) v_m^{(2g+2)},
\end{equation}
and then take one final step to reach
 \[ v_{m-1} = T_{i, j, m-1} + D^t_H((m-1)p)^{-1} M^t_H((m-1)p) v_m' \in W_{(m-1)p-1, t}. \]

If all of the above computations are performed to infinite precision, then it follows from Propositions \ref{prop:single-horizontal-reduction} and \ref{prop:horizontal-reduction} that if $v_m$ satisfies \eqref{eq:horizontal-pile}, then also $v_{m-1}$ also satisfies \eqref{eq:horizontal-pile}, and then by induction also $v_0$ satisfies \eqref{eq:horizontal-pile}.

Now we analyse the propagation of errors. To facilitate the analysis, we introduce the following terminology. Suppose that $v$ is a vector of length $2g+1$, with coordinates in $R_1$. Let $\error{v}$ denote the error term associated to $v$. That is, $\error{v}$ is the difference between the value for $v$ stored by the machine and the value that would have been obtained for $v$ if all computations had been performed to infinite precision. We will say that $v$ is \emph{$\ell$-correct} if:
\begin{itemize}
\item the $\ell$-th coordinate of $v$ is divisible by $p$;
\item the $\ell$-th coordinate of $\error{v}$ is divisible by $p^{N+1}$; and
\item the remaining coordinates of $\error{v}$ are divisible by $p^N$.
\end{itemize}
Note that $v_{i+(2g+1)j+1}$ is $1$-correct, since its first coordinate is simply $B_{j, (2g+1)j}$, which has been computed in $R_1$ and is divisible by $p$, and the other coordinates are all zero. The following series of claims together show that if $v_m$ is $1$-correct, then also $v_{m-1}$ is $1$-correct. Consequently $v_0$ is $1$-correct, and in particular $w_{i, j}$ is computed correctly to precision $p^N$.

\begin{claim}
\label{claim:1}
Let $1 \leq \ell \leq 2g$. If $v_m^{(\ell)}$ is $\ell$-correct, then $v_m^{(\ell+1)}$ is $(\ell+1)$-correct.
\end{claim}
\begin{proof}
We first examine the form of the matrix $M^t_H(mp-\ell)$. Let $P(x)$, $C(s, x)$ and $C_h(s)$ be the polynomials introduced in Proposition \ref{prop:single-horizontal-reduction}. We are taking $s = mp-\ell \equiv -\ell \pmod p$ and $t \equiv -1/2 \pmod p$, so
 \[ C(x, s) \equiv -2\ell P(x) + 2xP'(x) \pmod p. \]
In particular the coefficient $C_\ell(s)$ of $x^\ell$ is zero modulo $p$, so the entry in the $(\ell+1)$-th row of the last column of $M^t_H(mp-\ell)$ is zero modulo $p$. Consequently the contribution to $v_m^{(\ell+1)}$ from the last entry of $v_m^{(\ell)}$ satisfies the required conditions.

Furthermore, it is clear from Proposition \ref{prop:single-horizontal-reduction} that the only other possibly nonzero entry in the $(\ell+1)$-th row appears in the $\ell$-th column. Therefore $v_m^{(\ell+1)}$ also receives a contribution from the $\ell$-th entry of $v_m^{(\ell)}$, which by hypothesis already satisfies the required conditions.

Finally, the denominator
 \[ D^t_H(mp-\ell) = (2g+1)(2t-1) - 2s \equiv -2((2g+1) - \ell) \pmod p \]
is a unit, so dividing by it does not disturb $\ell$-correctness.
\end{proof}

\begin{claim} If $v_m^{(2g+1)}$ is $(2g+1)$-correct, then $v_m^{(2g+2)}$ is correct modulo $p^N$.
\end{claim}
\begin{proof}
Let $w = M_H^t(mp-2g-1)v_m^{(2g+2)}$. We have
\begin{equation}
\label{eq:blahblahblah}
 D_H^t(mp-2g-1) = (2g+1)(2t-1) - 2(mp - 2g - 1) \equiv 0 \pmod p,
\end{equation}
so by the definition of $M^t_H$, the first $2g$ columns of $M_H^t(mp-2g-1)$ are zero modulo $p$. Since the first $2g$ coordinates of $v_m^{(2g+1)}$ are correct modulo $p^N$, the contribution they make to $w$ is divisible by $p$ and correct modulo $p^{N+1}$. The contribution from the last coordinate of $v_m^{(2g+1)}$ is by hypothesis already divisible by $p$ and correct modulo $p^{N+1}$.

We deduce that $w$ is divisible by $p$ and correct modulo $p^{N+1}$. It suffices now to show that the valuation of $D_H^t(mp-2g-1)$ is precisely $1$. Since we know it is odd and divisible by $p$, we must bound its absolute value below $p^2$. From \eqref{eq:blahblahblah} and the definition of $t$ we have
 \[ D_H^t(mp-2g-1) = ((2g+1)(2j+1) - 2m)p, \]
and then the desired result follows from \eqref{eq:p-assumption}, since $0 \leq m \leq (2g+1)(j+1)-1$.
\end{proof}

\begin{claim}
If $v_m^{(2g+2)}$ is correct modulo $p^N$, then so is $v_m'$.
\end{claim}
\begin{proof}
By \eqref{eq:transport-to-v-prime} it suffices to show that $D^t_H((m-1)p, mp-2g-2)$ is a unit. The latter quantity is
 \[ \prod_{s=(k-1)p+1}^{kp-2g-2} (2g+1)(2t-1) -2s \equiv \prod_{s=1}^{p-2g-2} -2((2g+1) + s) \pmod p \]
since $t \equiv -1/2 \pmod p$, so it is a unit.
\end{proof}
\begin{rem}
In the above proof, we only needed the values of $M^t_H((m-1)p, mp-2g-2)$ and $D^t_H((m-1)p, mp-2g-2)$ modulo $p^N$, not modulo $p^{N+1}$. This is why it is possible to do almost all of the work in the horizontal reductions using only $N$ digits.
\end{rem}

\begin{claim}
If $v_m'$ is correct modulo $p^N$, then $v_{m-1}$ is $1$-correct.
\end{claim}
\begin{proof}
The same argument used in the proof of Claim \ref{claim:1} shows that the first row of $M_H^t((m-1)p)$ is entirely zero modulo $p$, and that $D_H^t((m-1)p)$ is a unit. Therefore the contribution to $v_{m-1}$ from $v_m'$ is $1$-correct. The contribution from $T_{i, j, m-1}$ is also $1$-correct.
\end{proof}

\subsection{Vertical reduction phase}

We first prove some lemmas that will be used to analyse the error propagation in the vertical reduction phase.

\begin{lem}
\label{lem:vertical-invertible}
If $t \equiv 1/2 \pmod p$, then $M_V(t)$ is invertible modulo $p$.
\end{lem}
\begin{proof}
Under the hypothesis on $t$, it follows from the definition of $M_V(t)$ that the entries of its $(i+1)$-th column are given by the coefficients of $S_i'(x)$. To show that $M_V(t)$ is invertible modulo $p$, it suffices to show that the $\reduce S_i'(x)$ are linearly independent over $\FF_q$. If $\sum_{i=0}^{2g-1} \lambda_i \reduce S_i'(x) = 0$ is some linear relation, then we may integrate (permissible, by \eqref{eq:p-assumption}) to obtain $\sum_{i=0}^{2g-1} \lambda_i \reduce S_i(x) = \lambda_{2g}$ for some $\lambda_{2g} \in \FF_q$. Multiplying this by $\reduce Q'(x)$, and using \eqref{eq:bezout}, we obtain
 \[ \sum_{i=0}^{2g-1} \lambda_i x^i \equiv \lambda_{2g} \reduce Q'(x) \pmod{\reduce Q(x)}. \]
But $1, x, \ldots, x^{2g-1}, \reduce Q'(x)$ are linearly independent in $\FF_q[x]/\reduce Q(x)$, since $Q'(x)$ has degree $2g$ and unit leading term (again due to \eqref{eq:p-assumption}). This forces every $\lambda_i = 0$.
\end{proof}

\begin{rem}
It would be interesting to characterise the values of $t$ for which $M_V(t)$ is singular modulo $p$. For instance, in the case of an elliptic curve, Example \ref{ex:elliptic} shows that $M_V(t)$ is singular precisely when $t \equiv 7/6 \pmod p$ or $t \equiv 5/6 \pmod p$. By studying the kernels and images of such maps, it may be possible to reduce the working precision in the vertical reduction steps from $p^{N+1}$ to $p^N$, as was done for the horizontal reductions.
\end{rem}

\begin{lem}
\label{lem:vertical-zero}
If $t_0 \equiv -1/2 \pmod p$, then $M_V(t_0, t_0 + p)$ is zero modulo $p$.
\end{lem}
\begin{proof}
Since $M_V(t_0, t_0 + p)$ modulo $p$ only depends on $t_0$ modulo $p$, we may assume that $t_0 = (p-1)/2$.

Let
 \[ X = D_V(t_0, t_0+p+1)^{-1} M_V(t_0, t_0 + p + 1) \]
be the reduction map from $W_{-1, t_0+p+1}$ to $W_{-1, t_0}$. First we will show that $pX$ is integral. It is easy to check that $p^2 X$ is integral, by inspecting the powers of $p$ dividing $D_V(t_0, t_0 + p + 1)$, but the integrality of $pX$ requires more work. The proof is very similar to the proof of \cite[Lemma 2]{kedlaya}. Let $\omega \in W_{-1, t_0+p+1}$, say
  \[ \omega = F(x) y^{-2(t_0+p+1)} dx/y, \]
where $F \in \ZZ_q[x]$ has degree at most $2g-1$. Let $\eta = X\omega$, and write
 \[ \eta = G(x) y^{-2t_0} dx/y \]
where $G \in \QQ_q[x]$ has degree at most $2g - 1$. We need to show that $p\eta$ is integral.

Since $X$ is a reduction map, $\omega$ and $\eta$ are cohomologous, and the discussion preceding Proposition \ref{prop:single-vertical-reduction} shows that $\omega - \eta = dH$ where
 \[ H = \sum_{t=t_0+1}^{t_0 +p+1} H_t(x) y^{-2t+1} \]
for some polynomials $H_t \in \QQ_q[x]$ of degree at most $2g$. We may now use the same argument as in the proof of \cite[Lemma 2]{kedlaya} (namely, comparing the $y$-expansions of $\omega$, $\eta$ and $dH$ around each root of $Q(x)$) to deduce that $m \eta$ is integral, provided that $m/(2t-1)$ is integral for $t_0 \leq t \leq t_0 + p + 1$. In particular $p\eta$ is integral, since we have assumed that $t_0 = (p-1)/2$.
 
Now we may finish the proof of the lemma. We have
 \[ X = D_V(t_0, t_0 + p + 1)^{-1} M_V(t_0, t_0 + p) M_V(t_0 + p + 1). \]
By Lemma \ref{lem:vertical-invertible} we know that $M_V(t_0 + p + 1)$ is invertible modulo $p$, so its inverse is integral. Rearranging, we obtain
 \[ M_V(t_0, t_0 + p) = D_V(t_0, t_0 + p + 1) X M_V(t_0 + p + 1)^{-1}. \]
Note that $D_V(t_0, t_0 + p + 1) = \prod_{t=t_0+1}^{t_0+p+1} (2t-1)$ is divisible by $p^2$, since the first and last factors in the product are zero modulo $p$. The integrality of $pX$ then implies that $M_V(t_0, t_0 + p)$ is zero modulo $p$.
\end{proof}

Now we may describe the vertical reduction phase. The input consists of the differentials $w_{i, j}$ computed via the horizontal reductions. The output will be a collection of differentials $w_i \in W_{-1, 0}$ for $0 \leq i < 2g$ that are cohomologous to $T_i$, and correct modulo $p^N$.

The first step is to compute the vertical reduction matrices
\[
  M_j = \begin{cases}
           M_V\left(0, \frac{p-1}2\right) & j = 0, \\
           M_V\left(\frac{(2j-1)p-1}2, \frac{(2j+1)p-1}2\right) & 1 \leq j < N,
     \end{cases}
\]
and similarly for $D_j$, using Theorem \ref{thm:bgs-generalisation}, with entries in $R_1$. The invertibility hypotheses of Theorem \ref{thm:bgs-generalisation} are satisfied, because the total reduction length $K$ satisfies
 \[ K = \frac{(2(N-1)+1)p - 1}2 < \frac{(2N-1)p}2. \]
The latter is bounded by $p^2/6$ from \eqref{eq:p-assumption}, so certainly $\sqrt K + 1 < p$.

For $j \geq 1$, observe that $D_j$ has valuation precisely $1$, because in the product
 \[ D_j = \prod_{t=\frac12((2j-1)p+1)}^{\frac12((2j+1)p-1)} (2t-1), \]
the only term divisible by $p$ is the first one, and \eqref{eq:p-assumption} implies that it is less than $p^2$. Furthermore, $M_j$ is zero modulo $p$ by Lemma \ref{lem:vertical-zero}. Since $M_j$ and $D_j$ have been computed modulo $p^{N+1}$, we can therefore compute the (integral) matrix
 \[ X_j = D_j^{-1} M_j \]
correctly modulo $p^N$. For the $j = 0$ case, the product for $D_0$ shows that it is a unit, so $X_0 = D_0^{-1} M_0$ may be computed modulo $p^N$ as well. Note that $X_j$ is the vertical reduction map from $W_{-1, \frac12((2j+1)p-1)}$ to $W_{-1, \frac12((2j-1)p-1)}$ for $j \geq 1$, and to $W_{-1, 0}$ for $j = 0$.

Now we fix $0 \leq i < 2g$, and show how to compute $w_i$. We define a sequence of differentials
\[
\begin{aligned}
 v_{N-1} & = w_{i, N-1} & & \in W_{-1, \frac12((2N-1)p-1)}, \\
 v_{N-2} & = w_{i, N-2} + X_{N-1} v_{N-1} & & \in W_{-1, \frac12((2N-3)p-1)}, \\
 \vdots \\
 v_0 & = w_{i, 0} + X_1 v_1 & & \in W_{-1, \frac12(p-1)}.
\end{aligned}
\]
Using Proposition \ref{prop:vertical-reduction}, one checks by induction that
 \[ v_m \sim \sum_{j \geq m} T_{i, j} \]
for each $1 \leq m \leq N-1$, and that the coefficients of $v_m$ are correct modulo $p^N$. Finally one puts $w_i = X_0 v_0 \in W_{-1, 0}$, which by Proposition \ref{prop:vertical-reduction} is cohomologous to $T_i$, and again its coefficients are correct modulo $p^N$.

\begin{rem}
In the case $N = 1$, it is only necessary to compute $M_0$ modulo $p$, rather than modulo $p^2$ as described above, since no divisions by $p$ are involved at all. It is not clear to the author whether a similar optimisation is available when $N > 1$.
\end{rem}

\subsection{Complexity analysis}
\label{sec:complexity-analysis}

\begin{proof}[Proof of Theorem \ref{thm:complexity}]

We first consider the time spent in the applications of Theorem \ref{thm:bgs-generalisation}, which will be the dominant step when $p$ is large compared to $N$ and $g$. For both $R_0 = \ZZ_{p^n}/(p^N)$ and $R_1 = \ZZ_{p^n}/(p^{N+1})$, basic ring operations (addition, multiplication) have bit-complexity $\OO(Nn \log p)$, and the costs of polynomial and matrix arithmetic over $R_i$ are given by $\polytime(d) = \OO(dNn\log p)$ and $\matrixtime(m) = \OO(m^\omega Nn\log p)$. For the horizontal reductions, for each of $N$ rows, we applied Theorem \ref{thm:bgs-generalisation} with $K = O(pN)$ and $m = O(g)$. Therefore each row costs $\OO(p^{1/2} N^{3/2} g^\omega n)$. For the vertical reductions, we applied Theorem \ref{thm:bgs-generalisation} once, also with $K = O(pN)$ and $m = O(g)$. Therefore the total time is $\OO(p^{1/2} N^{5/2} g^\omega n)$.

Now we estimate the time for the remaining steps, which for sufficiently large $p$ will be negligible.

Computing the coefficients $C_{j, r}$ in Proposition \ref{prop:frob-diff-approx} requires only $O(N^2 g^2)$ ring operations, even if naive polynomial multiplication is used. In the formulae for the $B_{j, r}$, computing all the necessary binomial coefficients requires $O(N^2)$ ring operations, and then computing all the $B_{j, r}$ requires $O(N^2 g)$ ring operations. Therefore computing the $B_{j, r}$ requires $O(N^2 g^2)$ ring operations altogether.

Solving \eqref{eq:bezout} for each $i$ requires $O(g^2)$ ring operations, even by the naive Euclidean extended GCD algorithm, so computing the coefficients of $M_V(t)$ needs at most $O(g^3)$ ring operations. Computing the coefficients of $M_H^t(s)$ for each of the $N$ required values of $t$ requires $O(Ng)$ ring operations.

In the horizontal reduction phase, computing the inverse of the Vandermonde matrix requires $O(N^3)$ ring operations. Then for each of $N$ rows we must perform the following steps. First, compute the values of $F^{(i)}(0)p^i/i!$, costing $O(N^2 g^2)$ ring operations. Then use these values to compute $M(k) = F(kp)$ for $O(Ng)$ values of $k$; for each $k$ this costs $O(N g^2)$ ring operations, so over all $k$ this costs $O(N^2 g^3)$. The total cost over all rows is $O(N^3 g^3)$ ring operations.

Finally we must account for the `single step' reductions during the horizontal reduction phase, as these were performed without the assistance of Theorem \ref{thm:bgs-generalisation}. Each matrix-vector multiplication requires only $O(g)$ ring operations, due to the sparsity of the matrices. For each of $N$ rows, for each of $O(Ng)$ values of $m$, and for each of $O(g)$ values of $i$, there are $O(g)$ such steps, for a total cost of $O(N^2 g^4)$ ring operations.

Altogether the cost is $O(N^3 g^4)$ ring operations, with corresponding bit-complexity $\OO(N^4 g^4 n \log p)$.
\end{proof}

\section{Sample computations}
\label{sec:samples}

The author implemented the main algorithm in C++, only for the case $n = 1$, using Victor Shoup's NTL library (\cite{shoup}, version 5.4) for the polynomial arithmetic. The implementation uses the middle product algorithm \cite{middle-product} for the key polynomial shifting step, as suggested in \cite[p.~1786]{bgs}; this was made trivial thanks to Shoup's wonderfully modular FFT code. The matrix multiplication steps use the naive $O(n^3)$ algorithm.

The source code is freely available under a GPL license from the author's web page, \texttt{http://math.harvard.edu/\textasciitilde{}dmharvey/}. The functionality has been made available in the SAGE computer algebra system (version 2.5.1 and later) \cite{sage}. An example session:
\begin{verbatim}
sage: R.<x> = PolynomialRing(ZZ)
sage: from sage.schemes.hyperelliptic_curves.frobenius import frobenius
sage: frobenius(p = 10007, N = 3, Q = x^5 + 2*x + 1)
 [844821791581 220205295882 761288372988 276316151941]
 [380371243619 656847071320 602083441024 781051879529]
 [435515877861 568305615656 204167847992  67069787872]
 [365277275232 293850471444 438804747301 298366229783]
\end{verbatim}

The following sample computations were performed on a 1.8 GHz AMD Opteron processor running Linux; many thanks to William Stein for offering this machine for the computations. The machine has 64GB of RAM and 16 cores, but only a single core was used. The compiler used was GCC 4.1.2 with optimisation flag \texttt{-O3}, and NTL was linked with the GMP library (version 4.2.1, with Pierrick Gaudry's AMD assembly patch) for the underlying integer arithmetic.

\subsection{Dependence on $p$}

Table \ref{tab:vary-p} shows the time used to compute the Frobenius matrix over a range of $p$ for the genus two curve $y^2 = x^5 - 11x^4 + 7x^3 - 5x^2 + 3x - 2$, with precision $N = 3$. From Theorem \ref{thm:complexity}, one expects the running time to approximately double for every four-fold increase in $p$.

\begin{table}[h]
\begin{tabular}{llllllll}
\hline
$p$            & time      & & $p$            & time      & & $p$            & time \\
\hline
$2^{14} + 27$  & 0.25 sec  & & $2^{24} + 43$  & 30.8 sec  & & $2^{34} + 23$  & 27.5 min   \\
$2^{16} + 1$   & 0.56 sec  & & $2^{26} + 15$  & 1.06 min  & & $2^{36} + 31$  & 1.00 hours   \\
$2^{18} + 3$   & 2.80 sec  & & $2^{28} + 3$   & 2.26 min  & & $2^{38} + 7$   & 2.62 hours   \\
$2^{20} + 7$   & 6.33 sec  & & $2^{30} + 3$   & 5.32 min  & & $2^{40} + 15$  & 6.39 hours   \\
$2^{22} + 15$  & 15.0 sec  & & $2^{32} + 15$  & 11.1 min  & & $2^{42} + 15$  & 13.5 hours  \\
\hline
\end{tabular}
\caption{Running times for $g = 2$ and $N = 3$}
\label{tab:vary-p}
\end{table}

\subsection{Near-cryptographic sizes} For the purposes of constructing secure cryptosystems, it is useful to be able to determine the zeta function of a hyperelliptic curve $C$ of low genus over a large prime field \cite[Ch.~23]{handbook}. In particular one hopes to find a curve whose Jacobian order $\#J(C/\FF_p)$ is prime (or is a prime multiplied by a very small integer) and sufficiently large.

For genus three and four, we ran our implementation on a single curve over the largest prime field that seemed feasible with the given hardware. We were able to determine the zeta function for a curve whose Jacobian approaches a cryptographically useful size, although there is still a gap to overcome. Handling a genus \emph{two} curve with a large enough Jacobian is clearly out of reach of this technique.

Thanks to Kiran Kedlaya for his assistance in using the MAGMA computer algebra system to perform some of the computations below.

\subsubsection{Genus three}

We computed the characteristic polynomial of Frobenius modulo $p$ for the curve
 \[ y^2 = x^7 + 17x^6 + 13x^5 + 11x^4 + 7x^3 + 5x^2 + 3x + 2 \]
defined over $\FF_p$ where $p = 2^{50} - 27$. The running time was 40 hours, and peak memory usage was 16 GB.

This determines $\#J(C/\FF_p)$ modulo $p$, within an interval of width $O(p^{3/2})$. The search space is only $O(p^{1/2})$, so MAGMA's baby-step/giant-step implementation is easily able to recover the Jacobian order. From this we inferred that the characteristic polynomial of Frobenius is
 \[ (X^6 + p^3) + a_1 (X^5 + p^2 X) + a_2 (X^4 + p X^2) + a_3 X^3, \]
where
\[
\begin{aligned}
 a_1 & = -8207566, \\
 a_2 & = 336549388766991, \\
 a_3 & = 17004180735172175425188. \\
\end{aligned}
\]
The order of the Jacobian over $\FF_p$ is
 \[ 1427247682301531613968301082755745957628851920 \sim 2^{150}. \]

\subsubsection{Genus four}

We computed the characteristic polynomial of Frobenius modulo $p^2$ for the curve
 \[ y^2 = x^9 - 23x^8 + 19x^7 - 17x^6 + 13x^5 - 11x^4 + 7x^3 - 5x^2 + 3x - 2 \]
defined over $\FF_p$ where $p = 2^{44} + 7$. The running time was 45 hours, and peak memory usage was 34 GB.

This does not pin down the zeta function precisely, but it produces a short list of four candidates, which we checked in MAGMA by testing which proposed Jacobian order $m$ satisfied $mP = 0$ for a number of random points $P$ defined over $\FF_p$. We found that the characteristic polynomial of Frobenius is
 \[ (X^8 + p^4) + a_1 (X^7 + p^3 X) + a_2 (X^6 + p^2 X^2) + a_3 (X^5 + p X^3) + a_4 X^4, \]
where
\[
\begin{aligned}
 a_1 & = 2394254, \\
 a_2 & = 29576915959850, \\
 a_3 & = 88182558522652238508, \\
 a_4 & = 536178748943545477971279916.
\end{aligned}
\]
The order of the Jacobian over $\FF_p$ is
 \[ 95780984339838343855809310281601230464609800042292722 \sim 2^{176}. \]

\bibliographystyle{amsalpha}
\bibliography{kedlaya-large-p}

\newcommand{\etalchar}[1]{$^{#1}$}
\providecommand{\bysame}{\leavevmode\hbox to3em{\hrulefill}\thinspace}
\providecommand{\MR}{\relax\ifhmode\unskip\space\fi MR }
\providecommand{\MRhref}[2]{%
  \href{http://www.ams.org/mathscinet-getitem?mr=#1}{#2}
}
\providecommand{\href}[2]{#2}
\begin{thebibliography}{CFA{\etalchar{+}}06}

\bibitem[BGS07]{bgs}
Alin Bostan, Pierrick Gaudry, and Eric Schost, \emph{Linear recurrences with
  polynomial coefficients and application to integer factorization and
  {C}artier--{M}anin operator}, SIAM Journal on Computing \textbf{36} (2007),
  no.~6, 1777--1806.

\bibitem[CC88]{chudnovsky}
D.~V. Chudnovsky and G.~V. Chudnovsky, \emph{Approximations and complex
  multiplication according to {R}amanujan}, Ramanujan revisited
  (Urbana-Champaign, Ill., 1987), Academic Press, Boston, MA, 1988,
  pp.~375--472.

\bibitem[CFA{\etalchar{+}}06]{handbook}
Henri Cohen, Gerhard Frey, Roberto Avanzi, Christophe Doche, Tanja Lange, Kim
  Nguyen, and Frederik Vercauteren (eds.), \emph{Handbook of elliptic and
  hyperelliptic curve cryptography}, Discrete Mathematics and its Applications
  (Boca Raton), Chapman \& Hall/CRC, Boca Raton, FL, 2006.

\bibitem[Edi03]{edixhoven}
Bas Edixhoven, \emph{Point counting after {K}edlaya}, EIDMA-Stieltjes Graduate
  course, Leiden (unpublished lecture notes),
  http://www.math.leidenuniv.nl/\textasciitilde{}edix/oww/mathofcrypt/carls\_e%
dixhoven/kedlaya.pdf (retrieved Oct 25th 2006), 2003.

\bibitem[GG03]{gaudry}
Pierrick Gaudry and Nicolas G{\"u}rel, \emph{Counting points in medium
  characteristic using {K}edlaya's algorithm}, Experiment. Math. \textbf{12}
  (2003), no.~4, 395--402.

\bibitem[HQZ04]{middle-product}
Guillaume Hanrot, Michel Quercia, and Paul Zimmermann, \emph{The middle product
  algorithm. {I}}, Appl. Algebra Engrg. Comm. Comput. \textbf{14} (2004),
  no.~6, 415--438.

\bibitem[Hub07]{hubrechts}
Hendrik Hubrechts, \emph{Quasi-quadratic elliptic curve point counting using
  rigid cohomology}, http://wis.kuleuven.be/algebra/hubrechts/ (retrieved May
  26th 2007), 2007.

\bibitem[Ked01]{kedlaya}
Kiran~S. Kedlaya, \emph{Counting points on hyperelliptic curves using
  {M}onsky-{W}ashnitzer cohomology}, J. Ramanujan Math. Soc. \textbf{16}
  (2001), no.~4, 323--338.

\bibitem[Ked04]{kedlaya2}
\bysame, \emph{Computing zeta functions via {$p$}-adic cohomology}, Algorithmic
  number theory, Lecture Notes in Comput. Sci., vol. 3076, Springer, Berlin,
  2004, pp.~1--17.

\bibitem[MST06]{mst}
B.~Mazur, W.~Stein, and J.~Tate, \emph{Computation of p-adic heights and log
  convergence}, Documenta Math. (Extra Volume: John H. Coates' Sixtieth
  Birthday) (2006), 577--614.

\bibitem[Sho07]{shoup}
Victor Shoup, \emph{{NTL}: A library for doing number theory},
  http://www.shoup.net/ntl/, 2007.

\bibitem[SJ05]{sage}
William Stein and David Joyner, \emph{Sage: System for algebra and geometry
  experimentation}, Communications in Computer Algebra (ACM SIGSAM Bulletin)
  \textbf{39} (2005), no.~2, 61--64.

\bibitem[Str69]{strassen}
Volker Strassen, \emph{Gaussian elimination is not optimal}, Numer. Math.
  \textbf{13} (1969), 354--356.

\end{thebibliography}

\end{document}